\documentclass{jnmp}
\usepackage{amsmath}

\JNMPnumberwithin{equation}{section}
\begin{document}
\renewcommand{\evenhead}{B~A Kupershmidt}
\renewcommand{\oddhead}{$q$-Analogs of Classical 6-Periodicity: From Euler to Chebyshev}
\thispagestyle{empty}

\FirstPageHead{10}{3}{2003}{\pageref{kupershmidt-firstpage}--\pageref{kupershmidt-lastpage}}{Article}

\copyrightnote{2003}{B~A Kupershmidt}

\Name{$\boldsymbol{q}$-Analogs of Classical 6-Periodicity: \\
From Euler to Chebyshev}
\label{kupershmidt-firstpage}

\Author{Boris A KUPERSHMIDT}

\Address{The University of Tennessee Space Institute, Tullahoma, TN  37388, USA\\
E-mail: bkupersh@utsi.edu}

\Date{Received November 19, 2002; Accepted January 03, 2003}

\begin{abstract}
\noindent
The sequence of period~6 starting with $1, 1, 0, -1, -1, 0$ appears
in many different disguises in mathematics.  Various $q$-versions of this sequence are
found, and their relations with Euler's pentagonal numbers theorem and Chebyshev polynomials
are discussed.
\end{abstract}

\begin{flushright}
\begin{minipage}{7.5cm}
\sf The motto on Cardinal Newman's tomb ought\\
to be the funeral motto of every Catholic,\\
{\it Ex umbris et imaginibus in veritatem},\\
Out of shadows and appearances into
the truth.\\[2mm]
\null \hfill Ronald Knox, {\it The Pastoral Sermons}
\end{minipage}
\end{flushright}

\section{Introduction}

The sequence of period 6, starting with
\begin{equation}
1, 1, 0, - 1, -1, 0,
\end{equation}
appears as the coefficients in the expansion
\begin{equation}
\frac{1}{1-x+x^2} = 1 +x - x^3 - x^4 + x^6 + x^7 - x^9 - x^{10} + x^{12}
+ x^{13} - \cdots.
\end{equation}
This is easy to see ([15, Ch.~V]):
\begin{gather}
\sum^\infty_{s=0} a_s x^s = \frac{1}{1-x + x^2} = \frac{1 + x}{1 +x^3}
= \frac{1}{1+x^3} + x \frac{1}{1+x^3}\nonumber\\
 \phantom{\sum^\infty_{s=0} a_s x^s}{}=
 \sum^\infty_{\ell =0} (-x^3)^\ell + x \sum^\infty_{\ell=0} (-x^3)^\ell,
\end{gather}
so that
\begin{subequations}
\begin{gather}
a_{s+3} = - a_s, \qquad s \in {\mathbb Z}_+, \\
a_0 = 1, \qquad a_1 = 1, \qquad a_2 = 0:
\end{gather}
\end{subequations}
the 6-periodicity results from 3-antiperiodicity.
\newpage

This is so far unremarkable.  However,
\begin{subequations}
\begin{gather}
\frac{1}{1-x+x^2} = \frac{1}{1-x(1-x)} = \sum^\infty_{m=0} x^m(1-x)^m \\
\phantom{\frac{1}{1-x+x^2}}{}= \sum_{m \geq 0} x^m \sum^m_{k=0}
 \left({m \atop k}\right) (-1)^k x^k =
\sum^\infty_{n=0} x^n \sum^{\lfloor n/2 \rfloor}_{k=0}
(-1)^k \left({n-k \atop k}\right).
\end{gather}
\end{subequations}
Thus, the binomial sums
\begin{subequations}
\begin{equation}
s_n = \sum^{\lfloor n/2 \rfloor}_{k=0} (-1)^k \left({n-k \atop k}\right)
\end{equation}
are equal to
\begin{equation}
1, 1, 0, -1, -1, 0 \qquad {\rm for} \qquad n \equiv 0, 1, 2, 3, 4, 5 \ \ ({\rm mod \ 6}).
\end{equation}
\end{subequations}

Such periodicity is very rare for binomial sums; see some examples and references
in~[12, 13].  The problem is:  can this periodicity be quantized, in the
sense of $q$-mathematics?

\section{Guessing an answer}

We need to find a suitable $q$-analog of the series $\sum\limits^\infty_{m=0}
x^m (1-x)^m$ (1.5a).  After some experimenting and rescaling, the
following sum suggests itself for consideration:
\begin{subequations}
\begin{gather}
S (x) = S(x;q) = \sum^\infty_{m=0} x^m (1 \;\dot{-} \; x)^m   \\
\phantom{S (x) = S(x;q)}{}= \sum^\infty_{m=0} x^m (x;q)_m,
\end{gather}
\end{subequations}
where
\begin{gather}
(a \; \dot{+}\; b)^m = \left[ \prod^m_{i =0} \left(a + q^i b\right)\right]/\left(a+q^m b\right) =
\left\{\begin{array}{ll}
(a+b)\cdots \left(a + q^{m-1} b\right), & m \in {\mathbb N},\vspace{1mm}\\
1, &  m=0
\end{array}\right.\!\!\!
\end{gather}
and (see [2, p.~487])
\begin{gather}
(x; q)_m = \left[\prod^m_{i \geq 0} \left(1 - q^i x\right) \right]/
\left(1-q^mx\right) =
 \left\{\begin{array}{ll}
(1 -x)\cdots \left(1 - q^{m-1} x\right), & m \in {\mathbb N},\vspace{1mm}\\
 1, & m = 0.
\end{array}\right.\!\!\!
\end{gather}
By Euler's formula [2, p. 492],
\begin{equation}
(1 \;\dot -\; x)^m = \sum^m_{k=0} \left[{m \atop k}\right] q
^{\left({k \atop 2}\right)} (-x)^k,
\end{equation}
where
\begin{equation}
\left[{m \atop k}\right] = \left[{m \atop k}\right]_q = \frac{[m] \cdots [m-k+1]}
{[1] \cdots [k]}, \qquad k \in {\mathbb N}; \qquad \left[{m \atop 0}\right] = 1,
\end{equation}
are the $q$-binomial coefficients, and
\begin{equation}
[m] = [m]_q = \left(1 - q^m\right)/(1-q)
\end{equation}
is a $q$-analog of the classical number (or object) $m$.

Substituting (2.4) into (2.1), we get
\begin{equation}
S(x) = \sum^\infty_{m=0} x^m (1\;\dot -\; x)^m = \sum^\infty
_{n=0} x^n \sum^{\lfloor n/2 \rfloor}_{k=0} (-1)^k
\left[{n-k \atop k}\right]_q q^{({k \atop 2})}.
\end{equation}
Thus, as a $q$-analog of the classical binomial sum (1.6a) we obtain
\begin{equation}
S_n = \sum^{\lfloor n/2 \rfloor}_{k=0} (-1)^k \left[{n-k \atop k}\right]
q^{({k \atop 2})}.
\end{equation}
Is it any good?  Calculating a few terms, we find:
\begin{equation}
S(x) = 1+x + 0 \cdot x^2 - qx^3 - q^2 x^4 + 0 \cdot x^5 + q^5 x^6 + q^7 x^7 + 0 \cdot x^8
+ O\left(x^9\right).
\end{equation}
This looks promissing, as the 3-antiperiodicity is preserved here, in $q$-clothes,
but there are not enough terms to guess the rule for the $q$-exponents.
We need some sort of a functional equation to determine that rule, assuming
it exists.

We proceed as follows.  Set
\begin{equation}
S(x, a) = S(x, a;q) = \sum^\infty_{m=0} a^m x^m (1 \; \dot -\;  x)^m.
\end{equation}
Since
\begin{equation}
(1 \;\dot -\; x)^{m+1} = (1 - x) (1 \; \dot -\; qx)^m,
\end{equation}
we find that
\begin{subequations}
\begin{equation}
S(x,a) = 1 + ax(1-x) \sum^\infty_{m=0} (ax)^m (1 \; \dot -\;  qx)^m =
1+ ax(1-x) S(qx, a/q).
\end{equation}
Therefore, if we write
\begin{equation}
S(x, a) = \sum^\infty_{n=0} S_n (a) x^n,
\end{equation}
\end{subequations}
the functional equation (2.12a) converts itself into
\begin{subequations}
\begin{gather}
S_{n+2} (a) = a S_{n+1} (a/q) q^{n+1} - a S_n(a/q)q^n, \qquad n \in {\mathbb Z}_+,
\\
S_1 (a) = a, \qquad S_0 (a) = 1.
\end{gather}
\end{subequations}
In particular, for $a =q$, we find:
\begin{subequations}
\begin{gather}
S_0 (q) = 1, \qquad  S_1 (q) = q, \\
S_{n+2} (q) = q^{n+1} (q S_{n+1} - S_n), \qquad
 S_n = S_n (1), \qquad n \in {\mathbb Z}_+.
\end{gather}
\end{subequations}
Using the expansion (2.9), we obtain then
\begin{gather}
S(x, q) = \sum^\infty_{m=0} (q x)^m (1 \; \dot -\; x)^m =
 1 + qx + \left(q^2 - q\right) x^2 - q^2 x^3 - q^5 x^4
\nonumber\\
\phantom{S(x, q) =} {}- \left(q^2 - q^7\right)x^5
+ q^7 x^6 + q^{12} x^7 + \left(q^{15} - q^{12}\right) x^8 + O \left(x^9\right).
\end{gather}
This is a bit less simple than (2.9) but still looks enticing.  Moreover,
comparing both expansions we can't fail to notice some coincidences and
regularities; the following general ansatz suggests itself:
\begin{subequations}
\begin{align}
& S_{6n} (q) = q^{x(n)},&& S_{6n} = q^{x(n)-2n}, &&&\\
& S_{6n+1} (q) = q^{y(n)}, && S_{6n+1} = q^{x(n)}, &&&\\
& S_{6n+2} (q) = q^{y(n)} \left(q^{2n+1} -1\right), && S_{6n+2} = 0, &&&\\
& S_{6n+3} (q) = - q^{u(n)}, && S_{6n+3} = -q^{u(n)-2n-1}, &&&\\
& S_{6n+4} (q) = -q^{v(n)}, && S_{6n+4} = -q^{u(n)}, &&&\\
& S_{6n+5} (q) = -q^{v(n)} \left(q^{2n+2} -1\right) , && S_{6n+5} =0. &&&
\end{align}
\end{subequations}
Making this ansatz compatible with the relations (2.14) leads to the
determination of all the unknown exponents $x(n)$,
$y(n)$, $u(n)$, $v(n)$  in formulae (2.16):
\begin{subequations}
\begin{gather}
y(n) = x(n) + 4n+1, \\
u(n) = x(n) + 6n+2, \\
v(n) = u(n) + 4n+3, \\
x(n+1) = 6n+5 +u (n)
\end{gather}
\end{subequations}
$\Rightarrow $
\begin{subequations}
\begin{gather}
x(n) = 6n^2 +n,  \\
y(n) = 6n^2 +5n+1,  \\
u(n) = 6n^2 +7n+2  \\
v(n) = 6n^2 +11n +5.
\end{gather}
\end{subequations}

These are our {\it two} -- conjectured so far -- candidates for a quantum
versions of the 6-periodic classical binomial sum (1.6).  Are there any other
candidates?  It seems unlikely, at least if we insist on having relatively
{\it compact} answers.  For example,
\begin{subequations}
\begin{gather}
S_4 (a) /a^2 = q - a [3] +a^2, \\
S_5 (a) /a^3 = q [3]-a [4] +a^2,
\end{gather}
\end{subequations}
and for $a = q^L$ we don't get anything attractive unless $L=0$ or 1.
\newpage

Let us now collect our conjectured formulae into a series form.  For
$a = 1$, we get
\begin{subequations}
\begin{gather}
S(x) = S (x,1) = \sum^\infty_{n=0} S_n x^n
\nonumber\\
{}=\sum^\infty_{n=0} \left(x^{6n} q^{6n^{2}-n} + x^{6n+1} q^{6n^{2}+n} - x
^{6n+3} q^{6n^2 + 5n+1} -  x^{6n+4} q^{6n^2 + 7n+2}\right)
\\
= \sum^\infty_{n=0} \left\{\left(x^{6n} q^{6n^{2}-n} - x^{6n+3} q^{6n^{2}+5n+1}
\right) + \left(x^{6n+1} q^{6n^{2}+n} - x^{6n+4} q^{6n^{2}+7n+2}\right)\right\}
\nonumber\\
= \sum^\infty_{n=0} (-1)^n \left\{x^{3n} q^{n(3n-1)/2} +
 x^{3n+1} q^{n(3n+1)/2} \right\}:
\\
\sum^\infty_{m=0} x^m (1 \;\dot{-}\; x)^m = \sum^\infty_{n=0}
(-x)^{3n} q^{n(3n-1)/2} \left(1 + xq^n\right).\tag{2.21}
\end{gather}
\end{subequations}
Similarly, for $a = q$ we obtain
\setcounter{equation}{21}
\begin{subequations}
\begin{gather}
S (x,q) = \sum^\infty_{n=0} S_n (q) x^n \nonumber\\
= \sum^\infty_{n=0} x^{6n} \left\{q^{6n^{2}+ n} +
 xq ^{6n^{2}+5n+1} + x^2 q^{6n^{2}+5n+1} \left(q^{2n+1} -1\right)\right.\nonumber\\
\left. {}-x^3 q^{6n^{2}+7n+2} - x^4 q^{6n^{2}+11n+5}
-x^5 q^{6n^{2}+11+5} \left(q^{2n+2} -1 \right) \right\}
\\
= \sum^\infty_{n=0} (-1)^n \left\{ x^{3n} q^{n(3n+1)/2} + x^{3n+1}
q^{n(3n+5)2} +x^{3n+2} q^{n(3n+5)/2} \left(q^{n+1}-1\right)\right\}:\!\!\!\!\!\!
\\
\sum^\infty_{m=0} q^m x^m (1 \;\dot -\; x)^m  = \sum^\infty_{n=0}
(-x)^{3n} q^{n(3n+1)/2} \left\{1+x q^{2n} + x^2 q^{2n} \left(q^{n+1} - 1\right)\right\}.\tag{2.23}
\end{gather}
\end{subequations}
In view of the functional equation (2.12), formula (2.21) implies formula (2.23).
In the next Section we shall prove the former formula.

\medskip

\noindent
{\bf Remark 2.24.} Formula (2.12a) shows that the parameter
\setcounter{equation}{24}
\begin{equation}
y = ax
\end{equation}
remains invariant under iteration.  Therefore, let us set
\begin{equation}
R(x, y) = R(x, y; q) = \sum^\infty_{m=0} y^m (1 \;\dot -\; x)^m.
\end{equation}
Then
\begin{equation}
R(x, y) = 1 + y (1-x) R(qx, y).
\end{equation}
Setting
\begin{gather}
R(x,y) = \sum^\infty_{\ell =0} c_\ell (y) x^\ell,
\end{gather}
we convert the functional equation (2.27) into
\setcounter{equation}{29}
\begin{gather}
c_0 = 1 + y c_0, \tag{{\rm 2.29a}}\\
c_{\ell+1} = yq^{\ell+1} c_{\ell+1} -yq^\ell c_\ell \qquad \Rightarrow \tag{{\rm 2.29b}}\\
c_\ell = \frac{(-y)^\ell q^{\left({\ell \atop 2}\right)}}{(1 \;\dot -\;y)^{\ell+1}},
\qquad \ell \in {\mathbb Z}_+ \qquad  \Rightarrow \\
R(x,y) = \sum^\infty_{m=0} y^m (1 \;\dot -\; x)^m  = \sum^\infty_{\ell =0}
\frac{(-xy)^\ell}{(1 \; \dot -\;  y)^{\ell+1}} q^{\left({\ell \atop 2}\right)} \qquad
\Rightarrow \\
S(x,a) = \sum^\infty_{m=0} (ax)^m (1 \;\dot -\; x)^m = \sum^\infty_{\ell =0}
\frac{(-ax^2)^\ell}{(1 \;\dot -\;  ax)^{\ell+1}} q^{\left({\ell \atop 2}\right)}.
\end{gather}
Neither of these expansions, however, is helpful for out task of proving formula (2.21).  We
need something completely different.  Notice that for $q = 1$, formula (2.31) yields:
\begin{gather}
\frac{ 1}{1 - y (1-x)} = \sum^\infty_{m=0} y^m (1-x)^m = \sum^\infty_{\ell =0}
\frac{(-xy)^\ell}{(1-y)^{\ell +1}} = \frac{1}{1-y}\frac{1}{\displaystyle 1+ \frac{xy}{1-y}}.
\end{gather}

\section{Euler's tower}

We have to prove formula (2.21).  In the form (2.20b), it is:
\begin{equation}
\sum^\infty_{m=0} x^m (1 \;\dot -\; x)^m  = \sum^\infty_{n=0}
(-1)^n \left\{x^{3n} q^{n (3n-1)/2} + x^{3n+1} q^{n(3n+1)/2}\right\}.
\end{equation}
This is an identity between formal power series in $x$, with coefficients
that are polynomials in~$q$.  We can not therefore let $x$ to be 1
(unlike $q$).
But let's get wild for a moment.  The RHS of formula (3.1) for $x=1$ is:
\begin{equation}
\sum^\infty_{n=0} (-1)^n \left\{q^{n(3n-1)/2} + q^{n(3n+1)/2} \right\} = 1+
\sum^\infty_{n=-\infty } (-1)^n q^{n(3n+1)/2} .
\end{equation}
This is essentially the RHS of Euler's famous identity, conjectured by him
in 1741:
\begin{gather}
\prod^\infty_{j=1} \left(1-q^j\right) = \sum^\infty_{n=-\infty} (-1)^n q^{n(3n+1)/2}.
\end{gather}
(See [15, Ch.~6] for an English translation of Euler's fascinating memoir.)
Nowadays, Euler's formula (3.3) is subsumed by the more general Jacobi triple
product identity
\begin{gather}
\prod^\infty_{n=1} \left(1 + Q^{2n-1} z\right) \left(1 + Q^{2n-1} z^{-1}\right) \left(1- Q^{2n}\right)
= \sum^\infty_{n=-\infty} Q^{n^{2}} z^n
\end{gather}
(see [16, p.~10 and p.~186]).  Euler's formula (3.3) results from (3.4) when
one makes the substitution
\begin{equation}
Q = q^{3/2}, \qquad z = - q^{1/2}.
\end{equation}

Neither Euler's formula (3.3) nor Jacobi's formula (3.4) seem of any help to
our identity~(3.1).  However, Euler's ingenious {\it Proof} of his
identity (see [24, p.~281]) has enough ingredients in it to establish
(3.1).  Weil describes Euler's Proof as ``another dazzling display of
algebraic virtuosity but quite elementary \dots''.

Euler sets
\begin{equation}
P_0 = \prod^\infty_{j=1} \left(1 - q^j\right)
\end{equation}
and notices that
\begin{gather}
P_0 = 1 - q - q^2 P_1,
\end{gather}
where
\begin{gather}
P_1 = \sum^\infty_{\nu = 0} q^\nu (1-q) \cdots\left (1-q^{\nu+1}\right).
\end{gather}
This follows from the easily proven by induction identity
\begin{gather}
\prod^n_{j=1} (1 - \alpha_j) = 1 - \alpha_1 - \sum^n_{k=2}
\alpha_k (1 - \alpha_1) \cdots (1 - \alpha_{k-1}), \qquad  n \geq 2,
\end{gather}
by setting $\alpha_j = q^j$, replacing $k$ by $\nu+2$, and then passing to the
limit $n \rightarrow \infty$.  Euler then introduces the infinite tower of
series $P_n$:
\begin{gather}
P_n = \sum^\infty_{\nu = 0} q^{\nu n} \left(1 - q^n\right) \cdots \left(1 - q^{n+\nu}\right), \qquad
n \in {\mathbb N},
\end{gather}
and shows by a {\it different argument} that
\begin{gather}
P_n = 1 - q^{2n+1} - q^{3n+2} P_{n+1}, \qquad n \in {\mathbb N}.
\end{gather}
Euler's argument for $n > 0$ is ingenious:
\begin{gather*}
P_n = \sum^\infty_{\nu = 0} q^{\nu n} \left(1 - q^n\right) \cdots \left(1 - q^{n+ \nu}\right) \\
= \left(1 - q^n\right) + q^n \left(1 - q^n\right) \left(1 - q^{n+1}\right) + \sum^\infty_{\mu = 2}
q^{\mu n} \left(1 - q^n\right) \cdots \left(1 - q^{n+\mu}\right)
\\
= 1 - q^n + q^n \left(1 - q^{n+1}\right) - q^{2n} \left(1 - q^{n+1}\right)
\\
{}+ \sum^\infty_{\mu=2} q^{\mu n} \left(1 - q^{n+1}\right) \cdots \left(1 - q^{n + \mu}\right) -
q^n \sum^\infty_{\mu = 2} q^{\mu n} \left(1 - q^{n+1}\right) \cdots \left(1 - q^{n+ \mu}\right)
\end{gather*}
\begin{gather*}
= 1 - q^{2n+1} - q^{2n} \left(1 - q^{n+1}\right) +
 \sum^\infty_{\nu=0} q^{\nu n+2n} \left(1 - q^{n+1}\right) \cdots \left(1 - q^{n + 2+\nu}\right)\\
{}- \sum^\infty_{\mu=1} q^{(\mu +1)n} \left(1 - q^{n+1}\right) \cdots \left(1 - q^{n + \mu}\right) +
q^{2n} \left(1 - q^{n+1}\right)\\
=1-q^{2n+1} + \sum^\infty_{\nu=0} q^{\nu (n+1)+2n-\nu}
\left(1 - q^{n+1}\right)\cdots \left(1 - q^{n + 1+ \nu}\right)  \left(1 - q^{n+2+\nu}\right)\\
{}- \sum^\infty_{\nu=0} q^{(\nu +2)n} \left(1 - q^{n+1}\right)\cdots\left(1 - q^{n + 1+\nu}\right)\\
= 1-q^{2n+1} + \sum^\infty_{\nu=0} q^{\nu(n+1)} \left(1 - q^{n+1}\right)\cdots
\left(1 - q^{n + 1 + \nu}\right)
\left\{q^{2n-\nu} \left(1 - q^{n+2+\nu}\right) - q^{2n- \nu}\right\} \\
= 1-q^{2n+1} \!-q^{3n+2} \sum^\infty_{\nu=0} q^{\nu(n+1)}\!
\left(1 - q^{n+1}\right) \cdots \left(1 - q^{n + 1+\nu}\right) = 1-
q^{2n+1} - q^{3n+2} P_{n+1}.\!
\end{gather*}
Collecting all the relations (3.11) together, we find
\begin{gather}
\prod^\infty_{j=1} \left(1 - q^j\right) = P_0 = 1-q-q^2 P_1 = 1-q - q^2 \left(1-q^3 - q^5
P_2\right)\nonumber\\
= 1-q - q^2 \left(1-q^3\right) + q^{2+5} \left(1 - q^5 - q^8 P_3\right) = \cdots\nonumber\\
= 1-q + \sum^\infty_{n=1} (-1)^n q^{2+5+\cdots  +(3n-1)} \left(1 - q^{2n+1} \right) \\
= 1-q + \sum^\infty_{n=1} (-1)^n q^{n(3n+1)/2} \left(1-q^{2n+1}\right) = \sum
^\infty_{n=0} (-1)^n q^{n(3n+1)/2} \left(1-q^{2n+1} \right) \\
= \sum^\infty_{n=0} (-1)^n q^{n(3n+1)/2} + \sum^\infty_{n=0} (-1)^{n+1}
q^{(3n^{2} + n + 4n+2)/2} = \sum^\infty_{n=-\infty} (-1)^n q^{n(3n+1)/2} .
\end{gather}
While we are at it, let's notice that a similar iteration yields
\begin{gather}
P_\ell = \sum^\infty_{n=0} (-1)^n q^{3\ell n+n(3n+1)/2} (1-q^{2\ell+2n+1}), \qquad
\ell \in {\mathbb Z}_+.
\end{gather}

From this, we can readily see that
\begin{gather}
\bar P_\ell : = 1 + q^\ell P_\ell = \sum^\infty_{n=0} (-1)^n \left\{q^{3 \ell n}
q^{n(3n-1)/2} + q^{(3n+1)\ell} q^{n(3n+1)/2} \right\}.
\end{gather}

It now only remains to compare $\bar P_\ell$ with
\setcounter{equation}{17}
\begin{subequations}
\begin{gather}
{\bar S}_\ell = S(q^\ell; q) = S (x; q) |_{x=q^{\ell}} ,  \qquad \forall \; \ell \in  {\mathbb N}:\tag{3.17}\\
{\bar S}_\ell = S (q^\ell; q) = \sum^\infty_{m=0} q^{\ell m} (1 \; \dot - \;  q^\ell) ^m =
 1 + \sum^\infty_{m=1} q^{\ell m} \left(1 - q^\ell\right) \cdots
\left(1 - q^{\ell + m - 1}\right),\!\!\!\!\\
\bar P_\ell = 1 + q^\ell \sum^\infty_{\nu=0} q^{\nu \ell} \left(1 - q^\ell\right) \cdots
\left(1 - q^{\ell +\nu}\right),
\end{gather}
\end{subequations}
and these are identical.  Thus,
\begin{equation}
\bar S_\ell = \bar P_\ell, \qquad  \ell \in {\mathbb N}.
\end{equation}
Comparing formulae (3.1) and (3.16), we get
\begin{gather}
\left.\left\{ \sum^\infty_{m=0} x^m (1 \;\dot -\;  x)^m \right\} \right|_{x=q^{\ell}} \nonumber\\
\qquad = \left.\left( \sum^\infty_{n=0} (-1)^n \left\{ x^{3n} q^{n (3n-1)/2} + x^{3n+1} q^{n(3n+1)/2}
\right\}\right)\right|_{x=q^{\ell}}, \qquad \forall \; \ell \in {\mathbb N}.
\end{gather}
This proves formula (3.1), because:

\medskip

\noindent
{\bf Lemma 3.21.} {\it Suppose $f \in {\mathbb C} [[x,y]]$ is formal power series in $x$, $y$,
and let $r(1)$, $r(2)$, $\ldots$ be an increasing sequence of positive intergers.  If
\setcounter{equation}{21}
\begin{equation}
f \left(y^{r(i)}, y\right) = f (x,y) \big|_{x=y^{r(i)}} = 0, \qquad \forall \; i \in {\mathbb N},
\end{equation}
then $f$ is identically zero.}

\medskip

\noindent
{\bf Proof.} If $f (y^{r(1)}, y)=0$ then $f$ is divisible by $x-y^{r(1)} $:
\begin{gather}
f(x,y) = \left(x-y^{r(1)}\right) f_1 (x,y).
\end{gather}
Continuing on, we find that
\begin{gather}
f(x,y) = \left\{\prod^k_{i=1} \left(x-y^{r(i)}\right)\right\} f_k (x,y)
\end{gather}
Thus, all homogeneous components of $f$ of degrees $< k$ vanish, and $k$ is arbitrary.
\hfill\qed

\section{Euler's tower revisited}

\begin{flushright}
\begin{minipage}{7cm}
\sf {\it Woman:} Your're a pugilist, arent't you?\\
{\it World champ boxer Rocky Graziano:}  Nah, \\
I'm  just a prizefighter.
\end{minipage}
\end{flushright}

The Euler relation (3.11)
\begin{equation}
P_n = 1-q^{2n+1} - q^{3n+2} P_{n+1}, \qquad n \in {\mathbb N},
\end{equation}
for the sequence of series
\begin{equation}
P_n = \sum^\infty_{\nu=0} q^{\nu n} \left(1-q^n\right) \cdots \left(1 - q^{n+\nu}\right),
\qquad n \in {\mathbb N},
\end{equation}
has an implicit gem hidden inside.  To make it explicit, notice that in the definition
of~$P_n$~(4.2) the index $n$ could be treated as a formal parameter and not necessarily as
a~positive integer.  So, set
\begin{equation}
q^n = x
\end{equation}
and define
\begin{equation}
{\mathcal P} (x) = {\mathcal P} (x;q) = \sum^\infty_{\nu=0} x^\nu (1-x)
\cdots \left(1-xq^\nu\right).
\end{equation}
Applying the Euler argument to ${\mathcal P}(x)$, we get:
\begin{gather}
{\mathcal P} (x) = \sum^\infty_{\nu=0} x^\nu (1-x) \cdots \left(1-xq^\nu\right)
= (1 -x) + x(1-x) (1 - xq) \nonumber\\
{}+\sum^\infty_{\mu=2} x^\mu (1-x) \cdots\left(1-xq^\mu\right)
= (1-x) + x (1-xq) - x^2 (1-xq) \nonumber\\
{}+ \sum^\infty_{\mu=2}  x^\mu (1-xq) \cdots \left(1 - xq^\mu\right)
- x \sum^\infty_{\mu=2} x^\mu (1-xq)
\cdots \left(1-xq^\mu\right) \nonumber\\
{}= 1-qx^2 - x^2 (1-xq) + \sum^\infty_{\nu=0} x^{\nu+2} (1 - xq) \cdots \left(1 - xq^{\nu+1}\right)
\left(1 - xq^{\nu+2} \right)\nonumber\\
{}- \sum^\infty_{\mu=1} x^{\mu+1} (1-xq) \cdots \left(1 - xq^\mu\right) + x^2 (1-xq)
\nonumber\\
{}= 1-qx^2 + x^2 \sum^\infty_{\nu=0} x^\nu (1-xq) \cdots
\left(1-xq^{\nu+1}\right) \left\{\left(1-xq^{\nu+2}\right) - 1\right \}\nonumber\\
{}= 1-qx^2 - x^3 q^2 \sum^\infty_{\nu=0} (xq)^\nu (1-xq) \cdots (1-xq \cdot q)^\nu =
1-qx^2 - q^2 x^3 {\mathcal P} (xq): \\
{\mathcal P} (x;q) = 1-qx^2 - q^2 x^3 {\mathcal P} (xq; q).
\end{gather}
Writing
\begin{gather}
{\mathcal P} (x;q) = \sum^\infty_{\ell = 0} c_\ell (q) x^\ell,
\end{gather}
we can transform the relation (4.6) into
\begin{subequations}
\begin{gather}
c_{\ell+3} = - q^{\ell + 2} c_\ell, \\
c_0 =1, \qquad c_1 = 0, \qquad c_2 = -q \qquad \Rightarrow \\
c_{3\ell} = (-1)^\ell q^{\ell (3 \ell +1)/2}, \tag{\rm 4.9a}\\
c_{3\ell+1} = 0, \tag{\rm 4.9b}\\
c_{3\ell-1} = (-1)^\ell q^{\ell (3\ell-1)/2}  \qquad \Rightarrow \tag{\rm 4.9c}\\
{\mathcal P}(x) = \sum^\infty_{\nu=0} x^\nu (1-x) \cdots \left(1 - xq^\nu\right)\tag{\rm 4.10a}\\
{}\qquad= 1 + \sum^\infty_{\ell=1} (-1)^\ell \left\{ q^{\ell(3\ell-1)/2} x^{3 \ell-1} + q^{\ell (3\ell
+1)/2} x^{3 \ell} \right\}.\tag{\rm 4.10b}
\end{gather}
\end{subequations}

If we notice that
\setcounter{equation}{10}
\begin{equation}
S(x) = 1 + x {\mathcal P} (x),
\end{equation}
then from formula (4.10b) we obtain at once that
\begin{subequations}
\begin{gather}
S (x) = \sum^\infty_{m=0} x^m (1 \;\dot -\; x)^m\\
\qquad{}= \sum^\infty_{\ell=0} (-1)^\ell \left\{x^{3 \ell} q^{\ell (3 \ell -1)/2} + x^{3 \ell +1}
q^{\ell (3 \ell + 1)/2} \right\},
\end{gather}
\end{subequations}
which is our conjectured formula (2.20b).  The functional equation (4.6) in the $S$-language
becomes:
\begin{gather}
S(x) = 1 + x - qx^3 S(qx);
\end{gather}
this equation by itself implies formula (4.12b) at once.

\medskip

\noindent
{\bf Remark 4.14.} The identity (4.10) is essentially due to Euler,
but it was repeatedly rediscovered many times
since.  See [1, p.~282] for more on that.

\medskip

\noindent
{\bf Remark 4.15.} What is the orgin of number~6,
and are there interesting periodic sequences of period other than~6?
I don't yet have a definite answer to these questions but suspect
that our 6-periodic sequence is a camouflaged group
of units of the field ${\mathbb Q}[{\sqrt{-3}}]$;
 this suggests that there exists a 4-periodic sequence attached to the units
$\pm 1$, $\pm i$ of ${\mathbb Q}[{\sqrt{-1}}]$,
and if a number field $K \supset {\mathbb Q}$ has a finite group of units then these can be
organized into a periodic (vector) sequence.

\section{Recurrence relations of second order}

The space of 3-antiperiodic sequences is 3-dimensional.  It contains a special 2-dimensional
subspace consisting of 2$^{nd}$-order recurrent sequences $\{u_n\}$ satisfying the relation
\begin{equation}
u_{n+2} = u_{n+1} - u_n.
\end{equation}
Indeed, for $u_0 = a$, $u_1 = b,$ the sequence $u_n$ starts as
\begin{equation}
a, \ b, \ b-a, \ -a, \ -b, \ -b \ + a \ | \ a, \ b, \ \ldots
\end{equation}
(In the language of characteristic polynomials: $\lambda^3 +1$
is divisible by $\lambda^2 - \lambda +1$).

Our sequence of binomial sums $s_n$ (1.6a)
\begin{equation}
s_n = \sum^{\lfloor n/2 \rfloor}_{k=0} (-1)^k \left({n-k \atop k }\right)
\end{equation}
is of this special form, with $a = b = 1$.  This is true a posteriori, from
comparing formulae~(5.2) and (1.6b), but can be easily seen directly from
the definition of $s_n$ (5.3):
\begin{gather}
s_{n+2}- s_{n+1} = \sum_{k \geq 0} (-1)^k \left\{ \left({n+2-k \atop k}\right)
- \left( {n+1-k \atop k}\right) \right\} = \sum (-1)^k \left( {n + 1-k \atop
k-1} \right)\nonumber\\
\qquad = \sum (-1)^{k+1} \left({n - k \atop k} \right) = - s_n.
\end{gather}

In quantizing $s_n$ in \S~2, we went along the route of generating functions,
quantizing the latter for the sequence $s_n$:
\begin{equation}
\sum^\infty_{n=0} s_n x^n = \sum^\infty_{m=0} x^m (1-x)^m.
\end{equation}
Generating functions make a very useful and powerful device (see [25].)
However, when one attempts to quantize a classical object, property, relation,
etc., the generating function methods may unnecessarily restrict one's choice.
A cautionary tale follows.

Let's try to quantize the {\it calculation} (5.4).  Set
\begin{equation}
f_n = f_{n|L} = f_{n|L} (q) = \sum^{\lfloor n/2 \rfloor}_{k=0} \left(-q^L\right)^k
q^{\theta (n,k)} \left[ {n - k \atop k} \right]_q
\end{equation}
with some unknown function $\theta(n,k)$ to be specified later on.  We shall use
the formulae
\begin{subequations}
\begin{gather}
\left[ {\alpha + 1 \atop r} \right] = q^r \left[{\alpha \atop r}\right] +
\left[{\alpha \atop r - 1} \right]\\
\phantom{\left[ {\alpha + 1 \atop r} \right]}{} = \left[{\alpha \atop r}\right] + q^{\alpha + 1 - r} \left[{\alpha \atop
r-1}\right].
\end{gather}
\end{subequations}

First, by formula (5.7a)
\begin{gather}
f_{n+2|L} = \sum_{k} \left(-q^L\right)^k q^{\theta (n + 2, k)} \left[{n + 2 - k \atop
k} \right] \nonumber\\
= \sum_k \left(-q^L\right)^k q^{\theta (n+2, k)} \left\{q^k \left[{n + 1 - k \atop
k}\right] + \left[{n+1-k \atop k-1}\right]\right\}\nonumber\\
= \sum \left(-q^{L+1}\right)^k q^{\theta (n+2, k)} \left[{n+1-k \atop k}\right] - q^L
\sum \left(-q^L\right)^k q^{\theta (n+2,k+1)} \left[{n - k \atop k}\right].
\end{gather}
This is an opaque expression.  If we simplify it by having $L=0$, demanding
that
\begin{equation}
\theta (n+2, k) + k = \theta (n+1, k)
\end{equation}
modulo rescaling of $f_n$ by a function of $n$, and then requiring that
\begin{equation}
\{\theta (n+2, k+1) - \theta (n, k)\} \quad \mbox{\rm is \ $k$-independent},
\end{equation}
we find that
\begin{equation}
\theta (n, k) = k(k-n),
\end{equation}
and formula (5.8) yields:
\begin{equation}
f_{n+2} = f_{n+1} - q^{-(n+1)} f_n.
\end{equation}
Changing $q$ into $q^{-1}$, we get:
\begin{gather}
f_{n+2} = f_{n+1} - q^{n+1} f_n, \qquad f_0 = f_1 = 1,
\\
f_n = \sum^{\lfloor n/2 \rfloor}_{k=0} (-1)^k q^{k^{2}} \left[{n - k
\atop k}\right]_q
\end{gather}
because
\begin{subequations}
\begin{gather}
[n]!_{q-1} = [n]! q^{-({n \atop 2})}, \\
\left[{n \atop r}\right]_{q^{-1}}  = \left[{n \atop r}\right]_q q^{r(r-n)}.
\end{gather}
\end{subequations}
The sequence $f_n$ is:
\begin{subequations}
\begin{gather}
1, \ 1, \ 1-q, \ 1-q - q^2, \ 1-q -q^2 -q^3 + q^4, \ 1 - q -q^2 - q^3 + q^5 + q^6, \nonumber\\
1-q- q^2 - q^3 + 2q^6 + q^7 + q^8 + q^9, \ \ldots
\end{gather}
and it doesn't appear interesting; however, its limit $f_\infty$ certainly
is:  from formula (5.14) with $n=\infty$, we get
\begin{equation}
f_\infty = \sum^\infty_{k=0} \frac{(-1)^k q^{k^{2}}}
{(1 \; \dot - \; q)^k}.
\end{equation}
\end{subequations}

Secondly, if we use formula (5.7b) instead of (5.7a), we  find (with $L = 0$):
\begin{gather}
f_{n+2} = \sum (-1)^k q^{\theta (n+2, k)} \left\{ \left[{n + 1 - k \atop
k}\right] + q^{n+2-2k} \left[{n + 1 - k \atop k-1}\right] \right\} \nonumber\\
= \sum (-1)^k q^{\theta (n + 2, k)} \left[{n + 1 - k \atop k} \right] -
q^n \sum (-1)^k q^{\theta (n+2, k+1)-2k} \left[{n - k \atop k}\right].
\end{gather}
The first summand in (5.17) implies that $\theta (n, k)$ is $n$-independent,
and the second summand yields
\begin{equation}
\theta(n, k) = k^2.
\end{equation}
Thus,
\begin{gather}
f_{n+2} = f_{n+1} - q^{n+1} f_n, \qquad
 f_0 = f_1 = 1, \\
f_n = \sum^{\lfloor n/2 \rfloor}_{k=0} (-1)^k q^{k^{2}} \left[{n - k \atop
k}\right]_q,
\end{gather}
which is the same as formula (5.14).

If we set
\begin{equation}
F(x) = \sum^\infty_{n=0} f_n x^n,
\end{equation}
then formula (5.19) converts itself into
\begin{equation}
F(x) (1-x) = 1-qx^2 F (qx),
\end{equation}
whence, by repeated iteration,
\begin{equation}
F(x) = \sum^\infty_{k=0} \frac{\left(-x^2\right)^k q^{k^{2}}}
{(1\; \dot -\;  x)^{k+1}}.
\end{equation}

Our original sequence $S_n$ (2.8),
\begin{equation}
S_n = \sum^{\lfloor n/2 \rfloor}_{k=0} (-1)^k \left[{n-k \atop k}\right]
q^{\left({k \atop 2}\right)}
\end{equation}
satisfies, as can be directly checked, the following $q$-analog of the
relation (5.1):
\begin{equation}
S_n = q^{\lfloor n/3 \rfloor} S_{n-1} -q^{\lfloor 2n/3 \rfloor-1}
S_{n-2};
\end{equation}
the more general sequence $S_n (a)$ (2.12):
\begin{equation}
S_n (a) = \sum^{\lfloor n/2 \rfloor}_{k=0} (-1)^k \left[{n - k \atop k}
\right]q^{\left({k \atop 2}\right)} a^{n-k},
\end{equation}
for $a=q$ satisfies the easily verified relation
\begin{equation}
S_n (q) = q^{\lfloor (n+2)/3\rfloor} S_{n-1} (q) - q^{\lfloor (2n+1)/3
\rfloor} S_{n-2} (q).
\end{equation}

\section{Chebyshev polynomials and series}

The special 3-antiperiodic relation (5.1)
\begin{equation}
u_{n+2} = u_{n+1} - u_n
\end{equation}
is satisfied by the values of the Chebyshev polynomials $T_n(x) $ and
$U_n (x)$ at $x=1/2$.  Let us recall the salient facts.

Originally, Chebyshev polynomials appeared as the least deviant monic
polynomials of a fixed degree, $n$ say, on the interval $-1 \leq x \leq 1$.  For the norm
\[
\|f\| = \max\limits_{-1\leq x \leq 1} |f(x)|,
\]
the answer is $T_n (x)/2^{n-1}$,
\begin{equation}
T_n (x) = \cos (n \theta), \qquad  x = \cos \theta,
\end{equation}
being the Chebyshev polynomial of the 1$^{st}$ kind; for the norm
\[
\|f\| = \int^1_{-1} |f(x)|dx,
\]
the answer is $U_n (x)/2^n, $
\begin{equation}
U_n (x) = \frac{\sin ((n+1) \theta))}{\sin \theta} , \qquad  x = \cos \theta,
\end{equation}
being the Chebyshev polynomial of the 2$^{nd}$ kind.

The following table is extracted from [11, p.~287]:
\begin{gather}
T_0 (x) =1, \nonumber\\
T_1 (x) = x, \nonumber\\
T_2 (x) = 2x^2-1, \nonumber\\
T_3 (x) = 4x^3 - 3x, \nonumber\\
T_4 (x) = 8x^4 - 8x^2 +1, \nonumber\\
T_5 (x) = 16x^5 - 20x^3+ 5x, \nonumber\\
T_6 (x) = 32 x^6 - 48^4 + 18x^2 -1, \nonumber\\
T_7 (x) = 64x^7 - 112 x^5 + 56x^3 - 7x, \nonumber\\
T_8 (x) = 128x^8 - 256x^6 + 160x^4 - 32x^2 + 1, \nonumber\\
T_9 (x) = 256x^9 - 576x^7 + 232x^5 - 120^3 +9x; \\
U_0 (x) = 1, \nonumber\\
U_1 (x) = 2x, \nonumber\\
U_2 (x) = 4x^2 -1, \nonumber\\
U_3 (x) = 8x^3 - 4x, \nonumber\\
U_4 (x) = 16x^4 - 12x^2 +1, \nonumber\\
U_5 (x) = 32x^5 - 32x^3 + 6x, \nonumber\\
U_6 (x) = 64x^6 - 80x^4 + 24x^2 -1, \nonumber\\
U_7 (x) = 128x^7 - 192x^5 + 80x^3 - 8x, \nonumber\\
U_8 (x) = 256x^8 - 448x^6 + 240x^4 - 40x^2 + 1, \nonumber\\
U_9 (x) = 512x^9 - 1024x^7 + 672x^5 - 106x^3 + 10x.
\end{gather}

The table (6.4) suggests that
\begin{equation}
T_p (x) \equiv x^p \quad ({\rm mod} \ p), \qquad p \ \mbox{\rm an odd  prime},
\end{equation}
which is true because of the easily verifiable formula
\begin{equation}
T_n (x) = \sum^{\lfloor n/2 \rfloor}_{k=0} \left({n \atop 2k}\right)
x^{n-2k} \left(x^2 - 1\right)^k,
\end{equation}
and this interesting congruence is listed as an exercise in  [2, p.~117];
there are more exercises on that page, about other interesting arithmetic
properties of the Chebyshev polynomials of the $1^{st}$ and $2^{nd}$ kind.
Some others properties not on that page are:
\begin{subequations}
\begin{gather}
T_{p^{2}} (x) \equiv T_p \left(x^p\right) \quad \left({\rm mod} \ p^2\right), \qquad p \
\mbox{\rm an odd prime}, \\
T_{2^{n+1}} (x) \equiv 1 \quad \left({\rm mod} \ 2^{2n+1}\right), \qquad n \in {\mathbb Z}_+,
\end{gather}
\end{subequations}
which follow from formula (6.6) and formula (6.17) below.  No doubt there
exist many more interesting arithmetic properties of the Chebyshev polynomials.

Such as.  Define the sequence $\{\gamma_n\}$ by the rule:
\begin{subequations}
\begin{gather}
\gamma_{n+2} = 2 \gamma_{n+1} + \gamma_n, \qquad \gamma_0 = 0, \qquad \gamma_1 = 1.
\end{gather}
Then
\begin{gather}
T_{2n} \left(\pm\sqrt{-1}\right) = 1 + (-1)^n \gamma^2_n, \qquad n \in {\mathbb Z}_+, \\
T_{2n+1} \left(\pm \sqrt{-1}\right) = \pm \sqrt{-1} \left\{-1+ (-1)^n
\left(\gamma_{n+1}^2 - \gamma^2_n\right)/2\right\},
\qquad n \in {\mathbb Z}_+.
\end{gather}
\end{subequations}

Most trigonometric formulae translate into formulae for the Chebyshev
polynomials in view of the representations (6.2) and (6.3).  For example,
formula
\[
\cos u + \cos v = 2 \cos \frac{u + v}{2} \cos \frac{u - v}{2},
\]
for $u = (n+2) \theta$, $v = n \theta; $ yields
\begin{subequations}
\begin{equation}
T_{n+2} (x) = 2x T_{n+1} (x) - T_n (x);
\end{equation}
formula
\[
\sin u + \sin v = 2 \sin \frac{u+v}{2} \cos \frac{u - v}{2}
\]
yields
\begin{equation}
U_{n+1} (x) = 2x U_n (x) - U_{n-1} (x),
\end{equation}
\end{subequations}
for $u = (n+2) \theta$, $v = n \theta;$ and
\begin{equation}
 U_n (x) - U_{n-2} (x) = 2T_n (x)
\end{equation}
for $u = (n+1) \theta$, $v = - (n-1) \theta$.  (Many more formulae involving Chebyshev
polynomials can be found in [3, \S~10.11], [17, 7], and [14, \S~5.7].)

Formulae (6.10) for $x = 1/2$ show that $\{T_n (1/2)\}$ and $\{U_n (1/2)\}$
are 3-antiperiodic sequences of the form (6.1), with
\begin{subequations}
\begin{gather}
T_0 (1/2) = 1, \qquad T_1 (1/2) = 1/2, \\
U_0 (1/2) = 1, \qquad U_1 (1/2) = 1.
\end{gather}
\end{subequations}
The preceding Sections thus have dealt with quantum aspects of Chebyshev
polynomials $U_n(x)$'s at just {\it one point} $x=1/2$.  We won't
attempt anything quantum in this Section, as this would be a rather
formidable undertaking that is better left to the interested reader.
Instead, we shall explore the problem of {\it interpolating} the function
$u(n) = u_n $ from integers to real and complex numbers.

We need a few more standard formulae.

Set temporatily ${\tilde U}_n (x) = (n+1)^{-1} dT_{n+1}/dx$.  Differentiating
formula (6.10a), we get:
\begin{gather}
(n+1) ({\tilde U}_{n+1} - 2x {\tilde U}_n + {\tilde U}_{n-1}) +
({\tilde U}_{n+1} -
{\tilde U}_{n-1} - 2 T_n) = 0.
\end{gather}
This equation defines ${\tilde U}_n$'s recursively.  From the table (6.4),
\begin{subequations}
\begin{gather}
{\tilde U}_0 = T_1^\prime =1=U_0, \\
{\tilde U}_1 = T_2^\prime/2 = 2x = U_1,
\end{gather}
\end{subequations}
and the $U_n$'s also satisfy the equation (6.13) by formulae (6.10b) and
(6.11).  Therefore, ${\tilde U}_n = U_n$:
\begin{gather}
U_n (x) = \frac{1}{n+1} \frac{d T_{n+1}}{dx}.
\end{gather}
(Alternatively, one gets formula (6.15) by differentiating formula
(6.2) with respect to $x$.)

Formulae (6.10) and (6.15) yield:
\begin{subequations}
\begin{gather}
T_n (1) = 1, \\
U_n (1) = n+1, \\
T_n ' (1) = n^2.
\end{gather}
\end{subequations}

Formula (6.2) implies that ([19, p. 45],  [21, 27])
\begin{gather}
T_n (T_m (x)) = T_{nm} (x), \qquad n, m \in {\mathbb N}.
\end{gather}
Thus, the polynomials $T_n$'s form a commutative semigroup; below we shall
see that this semigroup is in fact a discrete part of a commutative 1-dimensional
{\it group} considered as either a formal group or as a group of holomorphic
automorphisms of a neighborhood of the complex plane around $z=1$.  Formula
(6.15) implies that $\{T_n$'s$\}$ are a more fundamental object than
$\{U_n$'s$\}$.

Now, from the definition (6.2),
\begin{subequations}
\begin{gather}
\frac{d^2 T_n}{d \theta^2} = - n^2 T_n.
\end{gather}
Since
\begin{gather}
\frac{d}{d \theta} = - \sin \theta \frac{d}{dx} , \qquad
\frac{d^2}{d \theta^2} = \left(1-x^2\right)
\frac{d^2}{dx^2} - x \frac{d}{dx},
\end{gather}
\end{subequations}
we get the classical equation
\begin{gather}
\left(1-x^2\right) T''_n - x T'_n + n^2 T_n = 0.
\end{gather}
The $T_n$ is the unique regular solution of the equation
\begin{subequations}
\begin{gather}
\left(1-x^2\right)y'' - xy' + n^2 y = 0
\end{gather}
satisfying the initial condition
\begin{gather}
y(1) = 1.
\end{gather}
\end{subequations}
In the variable
\begin{gather}
z = \frac{1 - x}{2},
\end{gather}
the equation (6.20a) takes the hypergeometric form
\begin{gather}
z (1-z) \frac{d^2 y}{dz^2} + \left(\frac 12 -z\right) \frac{dy}{dz} + n^2 y=0,
\end{gather}
and we recover another classical result:
\begin{gather}
T_n (x) = F\left(-n, n, \frac{1}{2}; \frac{1 -x}{2}\right).
\end{gather}
(The point $x=\frac 12$ directs us to the continuous 3-antiperiodic family
$u (\alpha) = F\!\left(- \alpha, \alpha, \frac{1}{2}; \frac{1}{4}\right)$.)

Formulae (6.4), (6.5), (6.10) imply that
\begin{subequations}
\begin{gather}
T_n (0) = U_n (0) = 0, \qquad  n \ \ {\rm odd}; \\
T_{2n} (0) = U_{2n} (0) = (-1)^n.
\end{gather}
\end{subequations}
The differential equation (6.19) then yields
\begin{gather}
T_n (x) = \frac{n}{2} \sum^{\lfloor n/2 \rfloor}_{k=0} (-1)^k
\frac{(n - k-1)!}{k!(n-2k)!} (2x)^{n-2k}, \qquad n \in {\mathbb N},
\end{gather}
and then formula (6.15) yields
\begin{gather}
U_n (x) = \sum^{\lfloor n/2 \rfloor}_{k=0} (-1)^k \left({n - k \atop k}
\right) (2x)^{n-2k}, \qquad n \in {\mathbb Z}_+;
\end{gather}
another classical pair of very useful formulae.

We now extend the index $n$ in $T_n (x)$ by allowing $n$ to be an arbitrary
complex number or a formal parameter:  we set
\begin{gather}
T_\alpha (x) = \sum^\infty_{k=0} c_{\alpha|k} (x-1)^k, \qquad  c_{\alpha|0}
=1,
\end{gather}
and require $T_\alpha$ to satisfy the differential equation
\begin{gather}
\left(1-x^2\right) T''_\alpha - x T'_\alpha + \alpha^2 T_\alpha = 0,
\end{gather}
together with the boundary condition
\begin{gather}
T_\alpha (1) = 1.
\end{gather}
This is equivalent to the relations
\begin{gather}
c_{\alpha | k+1} = \frac{\alpha^2 - k^2}{(k+1) (2k+1)} c_{\alpha|k} ,
\qquad k \in {\mathbb Z}_+ \quad \Rightarrow\\
c_{\alpha | k+1} = \frac{1}{(k+1)! (2k+1)!!} \prod^k_{i=0} \left(\alpha^2 - i^2\right),
\qquad c_{\alpha|0}=1.
\end{gather}
Since
\begin{gather}
\lim_{k \rightarrow \infty} |c_{\alpha|k+1} / c_{\alpha | k} | = 1/2,
\end{gather}
the series $T_\alpha (x)$ (6.27) converges for
\begin{gather}
|x-1|<2.
\end{gather}
The series terminates only when $\alpha $ is an integer.  It is readily
verified that
\begin{gather}
T_{\alpha +2} - 2 ( 1 + (x-1)) T_{\alpha +1} + T_\alpha = 0, \qquad
\forall\; \alpha,
\end{gather}
and that
\begin{gather}
U_\alpha - U_{\alpha -2} = 2 T_\alpha, \qquad  \forall\; \alpha,
\end{gather}
where
\begin{gather}
U_\alpha (x) = \frac{1}{\alpha +1} \frac{dT_\alpha}{dx}, \qquad \forall \;
\alpha.
\end{gather}
From this it follows that
\begin{gather}
U_{\alpha +2} - 2 ( (1+ (x-1)) U_{\alpha +1} + U_\alpha =0, \qquad
\forall\; \alpha.
\end{gather}
Also, formula (6.30) shows that
\begin{gather}
T_{- \alpha} = T_\alpha, \qquad \forall \; \alpha.
\end{gather}

Finally,
\begin{subequations}
\begin{gather}
T_\alpha (T_\beta (x)) = T_{\alpha \beta} (x), \qquad \forall \; \alpha, \beta.
\end{gather}
Indeed, the difference
\begin{gather*}
T_\alpha (T_\beta (x)) - T_{\alpha \beta} (x),
\end{gather*}
for each fixed degree $k$ of $(x-1)^k$, is a {\it polynomial} in $\alpha$,
$\beta$ which vanishes, by (6.17), for all $\alpha, \beta \in {\mathbb N} \times
{\mathbb N}$.  Therefore, this polynomial is identically zero.

In particular, since $T_1(x) = x =1 + (x-1)$ plays the role of the
identity in the group $\{T_\alpha\}$, we have:
\begin{gather}
(T_\alpha)^{-1} = T_{1/\alpha}.
\end{gather}
\end{subequations}
Most of the classical formulae for Chebyshev polynomials survive the constructed
extension -- either by the argument given above to prove formula (6.39a), or by
observing that
\begin{gather}
T_\alpha (x) = F\left(-\alpha, \alpha, \frac{1}{2}, \frac{1 -x}{2}\right), \qquad
\forall \; \alpha.
\end{gather}

\noindent
{\bf Remark 6.41.} The theme of parameter extension is quite common
in mathematics.  The {\it methods} for constructing such extensions vary
wildly.  I chose in this Section the group property of $T_n$'s.  A more
direct approach is to extend the parameters' meaning in the hypergeometric
representations, as in formulae (6.23) and (6.40).  Then there are various {\it relations}
one can extend.  For example, the Legendre polynomials $P_n(x)$ enter into
the following relations with the $T_n (x)$'s ([14, p. 261]):
\setcounter{equation}{41}
\begin{subequations}
\begin{gather}
\left(n + \frac{1}{2}\right) (1 + x)^{1/2} \int^x_{-1} P_n (t) (x-t)^{-1/2} dt = T_n
(x) + T_{n-1} (x), \\
\left(n+\frac{1}{2}\right) (1-x)^{1/2} \int^1_x P_n (t) (t-x)^{-1/2} dt = T_n (x) - T_{n+1}
(x),
\end{gather}
\end{subequations}
where $|x|<1$, $n \in {\mathbb N}$, and the integrals are taken as Cauchy principal
values.

\medskip

\noindent
{\bf Remark 6.43.} Treatises and collections of formulae on special
functions should always be read with a pen, pencil, and salt shaker ready.
For example, the book [9] lists on p.~206 the following formula:
\setcounter{equation}{43}
\begin{gather}
U_{n+1} (x) = (1-x)^{1/2} \sum^{\lfloor n/2 \rfloor}_{k=0}
\frac{(-1)^k (n-k)!}{k! (n-2k)!} (2x)^{n-2k}.
\end{gather}
Misprints are widespread, common, and not unexpected.  Some errors are
more subtle, and propagate from one book to another for {\it many} years.
(See an example of 186 years old errors in~[20]).

\medskip

\noindent
{\bf Remark 6.45.} The Chebyshev polynomials $\{T_n (x)\}$ and the system $\{x^n\}$
are essentially the only two {\it polynomial} systems satisfying the semigroup/group
law (see [18, 4]).  More general systems, if any exist, could be found only in {\it series}.  For
example, if one restricts oneself to the hypergeometric (classical or basic) groups, one
looks for triples of functions
\setcounter{equation}{45}
\begin{subequations}
\begin{gather}
a(x) = {_p{\cal{F}}_q} (\ldots; \omega_1 (x-1)),\\
b(x) = {_r{\cal{F}}_s} (\ldots; \omega_2 (x-1)), \\
c(x) = {_P{\cal{F}}_Q} (\ldots; \Omega (x-1)),
\end{gather}
such that
\begin{gather}
a(b(x)) = c(x).
\end{gather}
\end{subequations}

\noindent
{\bf Remark 6.47.} Suppose one quantizes the Chebyshev polynomials $U_n (x)$ by
any of the standard methods of $q$-hypergeometric (= basic) series (see [2, 5, 6, 8].)  Will
then the sequences $S_n$ and $S_n(q)$ constructed in \S~2 appear as ``values at $x=1/2$''
of the quantized polynomials?  I don't see how, and think it is very unlikely, as one
glance at formulae~(5.25) and (5.27) will show.  There is a deep mystery hidden here.

\medskip

\noindent
{\bf Remark 6.48.}
 The proper meaning of the Chebyshev polynomials, as the generators of a~semigroup of
polynomical maps of a vector space, has been discovered by Veselov~[22,~23] and by
Hoffman and Withers [26, 10], with such semigroups
parametrized by affine Weyl groups.  It would be very interesting to find
the analogs of the one-dimensional point
$x = 1/2$ and periodic sequences, both classical and quantum, for this more general set-up.

\medskip

\begin{flushleft}
\begin{minipage}{8cm}
\sf I hope that posterity will judge me kindly, \\
not only as to the things which I have explained, \\
but also to those which I have intentionally \\
omitted so as to leave to others the pleasure\\
of discovery. \\[1mm]
\null\qquad\qquad\qquad   Descartes
\end{minipage}
\end{flushleft}

\subsection*{Acknowledgements}

I thank Alexander Veselov for very interesting comments
and the anonymous referee for useful suggestions.

\label{kupershmidt-lastpage}

\end{document}